%% file: main.tex
\documentclass[ijoc,nonblindrev]{informs3-empty} 

\OneAndAHalfSpacedXII 

\usepackage{listings}

\usepackage{natbib}
 \bibpunct[, ]{(}{)}{,}{a}{}{,}%
 \def\bibfont{\small}%
\TheoremsNumberedThrough     

\EquationsNumberedThrough    

\usepackage{algorithm}
\usepackage{algpseudocode}
\usepackage{booktabs}
\usepackage{tikz}
\usepackage{tikz-qtree}
\usepackage{multicol}
\usetikzlibrary{shapes}
\usepackage{fancyvrb}
\DefineVerbatimEnvironment%
 {MyCode}{Verbatim}
 {frame=single,baselinestretch=0.6,xleftmargin=0.5in,xrightmargin=0.4in,rulecolor=\color{gray},
 commandchars=\\\{\}}

\algrenewcommand\algorithmicforall{\textbf{For each}}
\usepackage{setspace}
\usepackage[hidelinks]{hyperref}
\begin{document}


\RUNAUTHOR{Lubin and Dunning}

\RUNTITLE{Computing in OR using Julia}

\TITLE{Computing in Operations Research using Julia}

\ARTICLEAUTHORS{%
\AUTHOR{Miles Lubin, Iain Dunning}
\AFF{MIT Operations Research Center, 77 Massachusetts Avenue, Cambridge, MA USA\\\EMAIL{mlubin@mit.edu}, \EMAIL{idunning@mit.edu}}
} 

\ABSTRACT{%
The state of numerical computing is currently characterized by a divide between 
highly efficient yet typically cumbersome low-level languages such as C, C++, and 
Fortran and highly expressive yet typically slow high-level languages such 
as Python and MATLAB. 
This paper explores how Julia, a modern programming language
for numerical computing which claims to bridge this divide 
by incorporating recent advances in language
and compiler design (such as just-in-time compilation), can be used
for implementing software and algorithms fundamental to the field of operations research, with a focus on mathematical optimization.
In particular, we demonstrate algebraic modeling for linear and nonlinear 
optimization
and a partial implementation of a practical simplex code.
Extensive cross-language benchmarks suggest that Julia is capable
of obtaining state-of-the-art performance.

}%



\maketitle

%

\input{intro}

\input{jump}
\input{nlp}
\input{simplex}


\ACKNOWLEDGMENT{%
This work would not be possible without the effort of the Julia team, Jeff Bezanson, Stefan Karpinski, Viral Shah, and Alan Edelman, as well as that of the larger community of Julia contributors. We acknowledge, in particular, Carlo Baldassi and Dahua Lin for significant contributions to the development of interfaces for linear programming solvers. We thank Juan Pablo Vielma for his comments on this manuscript which substantially improved its presentation. M. Lubin was supported by the DOE Computational Science Graduate Fellowship, which is provided under grant number DE-FG02-97ER25308. 
}

%
%
%


\bibliographystyle{ijocv081} 
\bibliography{refs} 


\end{document}

%% file: intro.tex
\section{Introduction}

Operations research and digital computing have grown hand-in-hand over the last 60 years, with historically large amounts of available computing power being dedicated to the solution of linear programs~\cite[]{bixby2002solving}. Linear programming is one of the key tools in the operations research toolbox and concerns the problem of selecting variable values to maximize a linear function subject to a set of linear constraints. This foundational problem, the algorithms to solve it, and its extensions form a large part of operations research-related computation. The purpose of this paper is to explore modern advances in programming languages that will affect how algorithms for operations research computation are implemented, and we will use linear and nonlinear programming as motivating cases.

The primary languages of high-performance computing have been Fortran, C, and C++ for a multitude of reasons, including their interoperability, their ability to compile to highly efficient machine code, and their sufficient level of abstraction over programming in an assembly language. These languages are compiled offline and have strict variable typing, allowing advanced optimizations of the code to be made by the compiler.

A second class of more modern languages has arisen that is also popular for scientific computing. These languages are typically interpreted languages that are highly expressive but do not match the speed of lower-level languages in most tasks. They make up for this by focusing on ``glue code'' that links together, or provides wrappers around, high-performance code written in C and Fortran. Examples of languages of this type would be Python (especially with the Numpy~\cite[]{numpy} package), R, and MATLAB. Besides being interpreted rather than statically compiled, these languages are slower for a variety of additional reasons, including the lack of strict variable typing.

Just-in-time (JIT) compilation has emerged as a way to have the expressiveness of modern scripting languages and the performance of lower-level languages such as C. JIT compilers attempt to compile at run-time by inferring information not explicitly stated by the programmer and use these inferences to optimize the machine code that is produced. Attempts to retrofit this functionality to the languages mentioned above has had mixed success due to issues with language design conflicting with the ability of the JIT compiler to make these inferences and problems with the compatibility of the JIT functionality with the wider package ecosystems.

Julia~\cite[]{JuliaArxiv} is a new programming language that is designed to address these issues. The language is designed from the ground-up to be both expressive and to enable the LLVM-based JIT compiler~\cite[]{lattner2004llvm} to generate efficient code. In benchmarks reported by its authors, Julia performed within a factor of two of C on a set of common basic tasks. The contributions of this paper are two-fold: firstly, we develop publicly available codes to demonstrate the technical features of Julia which greatly facilitate the implementation of optimization-related tools. Secondly, we will confirm that the aforementioned performance results hold for realistic problems of interest to the field of operations research.

This paper is not a tutorial. We encourage interested readers to view the language documentation at \href{http://julialang.org/}{\small\texttt{julialang.org}}. An introduction to Julia's syntax will not be provided, although the examples of code presented should be comprehensible to readers with a background in programming. The source code for all of the experiments in the paper is available in the online supplement\footnote{\url{http://www.mit.edu/~mlubin/juliasupplement.tar.gz}}. \textit{JuMP}, a library developed by the authors for mixed-integer algebraic modeling, is available directly through the Julia package manager, together with community-developed low-level interfaces to both Gurobi and the COIN-OR solvers Cbc and Clp for mixed-integer and linear optimization, respectively.

The rest of the paper is organized as follows. In Section~\ref{sec:JuMP}, we present the package \mbox{JuMP}. In Section~\ref{sec:NLP}, we explore nonlinear extensions.
In Section~\ref{sec:Simplex}, we evaluate the suitability of Julia for low-level implementation of numerical optimization algorithms by examining its performance on a realistic partial implementation of the simplex algorithm for linear programming.

%% file: jump.tex
\section{JuMP}
\label{sec:JuMP}

Algebraic Modeling Languages (AMLs) are an essential component in any operations researcher's toolbox. AMLs enable researchers and programmers to describe optimization models in a natural way, by which we mean that the description of the model in code resembles the mathematical statement of the model. AMLs are particular examples of domain-specific languages (DSLs) which are used throughout the fields of science and engineering. 

One of the most well-known AMLs is AMPL~\cite[]{AMPLBook}, a commercial tool that is both fast and expressive.
This speed comes at a cost: AMPL is not a fully-fledged modern programming language, which makes it a less than ideal choice for manipulating data to create the model, for working with the results of an optimization, and for linking optimization into a larger project.

Interpreted languages such as Python and MATLAB have become popular with researchers and practitioners alike due to their expressiveness, package ecosystems, and acceptable speeds.
Packages for these languages that add AML functionality such as YALMIP~\cite[]{YALMIP} for MATLAB and PuLP~\cite[]{PULP} and Pyomo~\cite[]{Pyomo} for Python address the general-purpose-computing issues of AMPL but sacrifice speed.
These AMLs take a nontrivial amount of time to build the sparse representation of the model in memory, which is especially noticeable if models are being rebuilt a large number of times, which arises in the development of models and in practice, e.g. in simulations of decision processes.
They achieve a similar ``look'' to AMPL by utilizing the \textit{operator overloading} functionality in their respective languages, which introduces significant overhead and inefficient memory usage. Interfaces in C++ based on operator overloading, such as those provided by the commercial solvers Gurobi and CPLEX, are often significantly faster than AMLs in interpreted languages, although they sacrifice ease of use and solver independence.

We propose a new AML, \textit{JuMP} (Julia for Mathematical Programming), implemented and released as a Julia package, that combines the speed of commercial products with the benefits of remaining within a fully-functional high-level modern language. 
We achieve this by using Julia's \textit{metaprogramming} features to turn natural mathematical expressions into sparse internal representations of the model without using operator overloading.
In this way we achieve performance comparable to AMPL and an order-of-magnitude faster than other embedded AMLs.

\subsection{Metaprogramming with macros}

Julia is a \textit{homoiconic} language: it can represent its own code as a data structure of the language itself. This feature is also found in languages such as Lisp. To make the concept of metaprogramming more clear, consider the following Julia code snippet:
\vspace{0.1in}
\begin{MyCode}[numbers=left]
macro m(ex)
    ex.args[1] = :(-)  # Replace operation with subtraction
    return esc(ex)     # Escape expression (see below)
end
x = 2; y = 5  # Initialize variables
2x + y^x      # Prints 29
@m(2x + y^x)  # Prints -21
\end{MyCode}

On lines 1-4 we define the \textit{macro} \texttt{m}. Macros are compile-time source transformation functions, similar in concept to the preprocessing features of C but operating at the syntactic level instead of performing textual substitution. When the macro is invoked on line 7 with the expression $2x + y^x$, the value of \texttt{ex} is a Julia object which contains a representation of the expression as a tree, which we can compactly express in Polish (prefix) notation as:
\[
(+, (\ast,2,x), (^\wedge,y,x))
\]
Line 2 replaces the $+$ in the above expression with $-$, where \texttt{:(-)} is Julia's syntax for the symbol $-$. Line 3 returns the \textit{escaped} output, indicating that the expression refers to variables in the surrounding scope.  Hence, the output of the macro is the expression $2x - y^x$, which is subsequently compiled and finally evaluated to the value $-21$.

Macros provide powerful functionality to efficiently perform arbitrary transformation of expressions. The complete language features of Julia are available within macros, unlike the limited syntax for macros in C. Additionally, macros are evaluated only once, at compile time, and so have no runtime overhead,  unlike \texttt{eval} functions in MATLAB and Python. (Note: with JIT compilation, ``compile time'' in fact occurs during the program's execution, e.g., the first time a function is called.) We will use macros as a basis for both linear and nonlinear modeling.  

\subsection{Language Design}\label{sec:mathprogdesign}

While we did not set out to design a full modeling language with the wide variety of options as AMPL, we have sufficient functionality to model any linear optimization problem with a very similar number of lines of code.
Consider the following simple AMPL model of a ``knapsack'' problem (we will assume the data are provided before the following lines):

\vspace{0.1in}
\begin{MyCode}
var x\{j in 1..N\} >= 0.0, <= 1.0;

maximize Obj:
  sum \{j in 1..N\} profit[j] * x[j];
		
subject to CapacityCon:
  sum \{j in 1..N\} weight[j] * x[j] <= capacity;
\end{MyCode}

The previous model would be written in Julia using JuMP with the following code:

\vspace{0.1in}
\begin{MyCode}
m = Model(:Max)

@defVar(m, 0 <= x[1:N] <= 1)

@setObjective(m, sum\{ profit[j] * x[j], j = 1:N \})

@addConstraint(m, sum\{ weight[j] * x[j], j = 1:N \} <= capacity)
\end{MyCode}

The syntax is mostly self-explanatory and is not the focus of this paper, but we draw attention to the similarities between the syntax of our Julia AML and existing AMLs. In particular, macros permit us to define new syntax such as \texttt{sum\{\}}, which is not part of the Julia language. 

\subsection{Building Expressions}

The model is stored internally as a set of rows until it is completely specified. 
Each row is defined by two arrays: the first array is the indices of the columns that appear in this row and the second contains the corresponding coefficients.
This representation is essentially the best possible while the model is being built and can be converted to a sparse column-wise format with relative efficiency.
The challenge then is to convert the user's statement of the problem into this sparse row representation as quickly as possible, while not requiring the user to express rows in a way that loses the readability that is expected from AMLs.

AMLs like PuLP achieve this with operator overloading.
By defining new types to represent variables, new definitions are provided for the basic mathematical operators when one or both the operands is a variable. 
The expression is then built by combining subexpressions together until the full expression is obtained.
This typically leads to an excessive number of intermediate memory allocations.
One of the advantages of AMPL is that, as a purpose-built tool, it has the ability to statically analyze the expression to determine the storage required for its final representation.
One way it may achieve this is by doing an initial pass to determine the size of the arrays to allocate, and then a second pass to store the correct coefficients.
Our goal with Julia was to use the metaprogramming features to achieve a similar effect and bypass the need for operator overloading.

\subsection{Metaprogramming implementation}

Our solution is similar to what is possible with AMPL and does not rely on operator overloading at all.
Consider the knapsack constraint provided in the example above. 
We will change the constraint into an equality constraint by adding a slack variable to make the expression more complex than a single sum. The \texttt{addConstraint} macro converts the expression
\vspace{0.1in}
\begin{MyCode}
@addConstraint(m, sum\{weight[j]*x[j], j=1:N\} + s == capacity)
\end{MyCode}
into the following code, transparently to the user:

\vspace{0.1in}
\begin{MyCode}
aff = AffExpr()

sumlen = length(1:N)
sizehint!(aff.vars, sumlen)
sizehint!(aff.coeffs, sumlen)
for i = 1:N
    addToExpression(aff, 1.0*weight[i], x[i])
end

addToExpression(aff, 1.0, s)

addToExpression(aff, -1.0, capacity)

addConstraint(m, Constraint(aff,"=="))
\end{MyCode}

The macro breaks the expression into parts and then stitches them back together as in our desired data structure.
\texttt{AffExpr} represents the custom type that contains the variable indices (\texttt{vars}) and coefficients (\texttt{coeffs}). 
In the first segment of code the macro pulls the indexing scheme from out of the sum and determines how long an array is required. 
Sufficient space to store the sum is reserved in one pass using the built-in function \texttt{sizehint!} before \texttt{addToExpression} (defined elsewhere) is used to fill it out. 
We use \textit{multiple dispatch} to let Julia decide what type of object \texttt{x[i]} is, either a constant or a variable placeholder, using its efficient built-in type inference mechanism. 
After the sum is handled, the single slack variable \textit{s} is appended and finally the right-hand-side of the constraint is set.
Note the invocation of \texttt{addToExpression} with different argument types in the last usage - this time two constants instead of a constant and a variable.
The last step is to construct the \texttt{Constraint} object that is essentially a wrapper around the expression and the sense. The function \texttt{addConstraint} is defined separately from the macro with the same name.
We note that our implementation is not as efficient as AMPL's can be; space for the coefficients of single varibles like \textit{s} is not preallocated, and so additional memory allocations are required; however, we still avoid the creation of many small temporary objects that would be produced with operator overloading.

\subsection{Benchmarks}

Different languages produce the final internal representation of the problem at different stages, making pure in-memory ``model construction'' time difficult to isolate.
Our approach was to force all the AMLs to output the resulting model in the LP and/or MPS file formats and record the total time from executing the script until the file is output.
We evaluated the performance of Julia relative to other AMLs by implementing two models whose size can be controlled by varying a parameter. Experiments were performed on a Linux system with an Intel Xeon E5-2650 processor.

\begin{enumerate}
	\item P-median: this model was used by~\cite{Pyomo} to compare Pyomo with AMPL. The model determines the location of $M$ facilities over $L$ possible locations to minimize the distance between each of $N$ customers and the closest facility. $C_i$ is a vector of customer locations that we generate randomly. In our benchmarks we fixed $M=100$ and $N=100$, and varied $L$. The results are in Table~\ref{tab:pmediantotal}, and show that JuMP is safely within a factor of two of the speed AMPL, comparable in speed to, if not occasionally faster than, Gurobi's C++ modeling interface, and an order of magnitude faster than the Python-based modeling languages. Note that JuMP does not need to transpose the naturally row-wise data when outputting in LP format, which explains the observed difference in execution times. In JuMP and AMPL, model construction was observed to consume 20\% to 50\% of the total time for this model.
\begin{alignat*}{1}
\min & \sum_{i=1}^{N}\sum_{j=1}^{L}\left|C_{i}-j\right|x_{ij}\\
\text{s.t. } & x_{ij}\leq y_{j}\quad i=1,\ldots,N,\quad,j=1,\ldots,L\\
 & \sum_{j=1}^{L}x_{ij}=1\quad i=1,\ldots,N\\
 & \sum_{j=1}^{L}y_{j}=M
\end{alignat*}

	\item Linear-Quadratic control problem (\textit{cont5\_2\_1}): this quadratic programming model is part of the collection maintained by Hans Mittleman~\cite[]{mittelmann2009benchmarks}. Not all the compared AMLs support quadratic objectives, and the quadratic objective sections of the file format specifications are ill-defined, so the objective was dropped and set to zero. The results in Table~\ref{tab:cont5total} mirror the general pattern of results observed in the p-median model. 
\begin{alignat*}{1}
\min_{y_{i,j},u_{i}} & \dots \\
\text{s.t. } & \frac{y_{i+1,j}-y_{i,j}}{\Delta t}=\frac{1}{2\left(\Delta x\right)^{2}}\left(y_{i,j-1}-2y_{i,j}+y_{i,j+1}+y_{i+1,j-1}-2y_{i+1,j}+y_{i+1,j+1}\right)\\
 & \qquad\qquad\qquad i=0,\dots M-1,\ j=1,\dots,N-1\\
 & y_{i,2}-4y_{i,1}+3y_{i,0}=0\qquad\quad i=1,\dots,M\\
 & y_{i,n-2}-4y_{i,n-1}+3y_{i,n}= (2\Delta x) (u_{i}-y_{i,n})\qquad\quad i=1,\dots,M\\
 & -1\leq u_{i}\leq1\qquad\quad i=1,\dots,M\\
 & y_{0,j}=0\qquad\quad j=0,\dots,N\\
 & 0\leq y_{i,j}\leq1\qquad\quad i=1,\dots,m,\ j=0,\dots,N\\
\text{where } & g_{j}=\frac{1}{2}\left(1-\left(j\Delta x\right)^{2}\right)
\end{alignat*}

\end{enumerate}

\newcolumntype{x}[1]{%
>{\raggedleft\hspace{0pt}}p{#1}}%

\begin{table}[ht]
	\centering
	\caption{P-median benchmark results. L is the number of locations. Total time (in seconds) to process the model definition and produce the output file in LP and MPS formats (as available).}
	\label{tab:pmediantotal}
	\begin{tabular}{lx{1.1cm}x{1.1cm}x{1.1cm}x{1.1cm}x{1.1cm}x{1.1cm}x{1.1cm}x{1.1cm}}
		\toprule
		&\multicolumn{2}{c}{JuMP/Julia}&\multicolumn{1}{c}{AMPL}&\multicolumn{2}{c}{Gurobi/C++}&\multicolumn{2}{c}{Pulp/PyPy}&\multicolumn{1}{c}{Pyomo}\tabularnewline
		\cmidrule(l){2-3}\cmidrule(l){4-4}\cmidrule(l){5-6}\cmidrule(l){7-8}\cmidrule(l){9-9}
		L  & LP & MPS & MPS & LP & MPS & LP & MPS & LP\tabularnewline
		\cmidrule(l){1-1}\cmidrule(l){2-2}\cmidrule(l){3-3}\cmidrule(l){4-4}\cmidrule(l){5-5}\cmidrule(l){6-6}\cmidrule(l){7-7}\cmidrule(l){8-8}\cmidrule(l){9-9}
		1,000 & 0.5 & 1.0 & 0.7 & 0.8 & 0.8 & 5.5 & 4.8 & 10.7\tabularnewline
		5,000 & 2.3 & 4.2 & 3.3 & 3.6 & 3.9 & 26.4 & 23.2 & 54.6\tabularnewline
		10,000 & 5.0 & 8.9 & 6.7 & 7.3 & 8.3 & 53.5 & 50.6 & 110.0\tabularnewline
		50,000 & 27.9 & 48.3 & 35.0 & 37.2 & 39.3 & 224.1 & 225.6 & 583.7\tabularnewline 
		\bottomrule\tabularnewline
	\end{tabular}

\end{table}

\begin{table}[ht]
	\centering
	\caption{Linear-quadratic control benchmark results. N=M is the grid size. Total time (in seconds) to process the model definition and produce the output file in LP and MPS formats (as available).}
	\label{tab:cont5total}
	\begin{tabular}{lx{1.1cm}x{1.1cm}x{1.1cm}x{1.1cm}x{1.1cm}x{1.1cm}x{1.1cm}x{1.1cm}}
		\toprule
		&\multicolumn{2}{c}{JuMP/Julia}&\multicolumn{1}{c}{AMPL}&\multicolumn{2}{c}{Gurobi/C++}&\multicolumn{2}{c}{Pulp/PyPy}&\multicolumn{1}{c}{Pyomo}\tabularnewline
		\cmidrule(l){2-3}\cmidrule(l){4-4}\cmidrule(l){5-6}\cmidrule(l){7-8}\cmidrule(l){9-9}
		N  & LP & MPS & MPS & LP & MPS & LP & MPS & LP\tabularnewline
		\cmidrule(l){1-1}\cmidrule(l){2-2}\cmidrule(l){3-3}\cmidrule(l){4-4}\cmidrule(l){5-5}\cmidrule(l){6-6}\cmidrule(l){7-7}\cmidrule(l){8-8}\cmidrule(l){9-9}
		250 & 0.5 & 0.9 & 0.8 & 1.2 & 1.1 & 8.3 & 7.2 & 13.3\tabularnewline
		500 & 2.0 & 3.6 & 3.0 & 4.5 & 4.4 & 27.6 & 24.4 & 53.4\tabularnewline
		750 & 5.0 & 8.4 & 6.7 & 10.2 & 10.1 & 61.0 & 54.5 & 121.0\tabularnewline
		1,000 & 9.2 & 15.5 & 11.6 & 17.6 & 17.3 & 108.2 & 97.5 & 214.7\tabularnewline
		\bottomrule\tabularnewline
	\end{tabular}
\end{table}

\subsection{Availability}
JuMP (\href{https://github.com/IainNZ/JuMP.jl}{\small\texttt{https://github.com/IainNZ/JuMP.jl}}) has been released with documentation as a Julia package. It remains under active development. We do not presently recommend its use in production environments. It currently interfaces with both Gurobi and the COIN-OR solvers Cbc and Clp for mixed-integer and linear optimization. Linux, OS X, and Windows platforms are supported.

%% file: nlp.tex
\pagebreak 
\section{Nonlinear Modeling}
\label{sec:NLP}

Commercial AMLs such as AMPL and GAMS~\cite[]{GAMS} are widely used for specifying large-scale nonlinear optimization problems, that is, problems of the form
\begin{align}
	\min_x\quad&f(x)\notag\\
	\text{subject to}\quad& g_i(x) \le 0\quad i = 1,\ldots,m,\label{eq:constr}
\end{align}
where $f$ and $g_i$ are given by closed-form expressions. 

Similar to the case of modeling linear optimization problems, open-source AMLs exist and provide comparable, if not superior, functionality; however, they may be significantly slower to build the model, even impractically slower on some large-scale problems. The user guide of CVX, an award-winning open-source AML for convex optimization built on MATLAB, states that ``CVX is \textit{not} meant for very large problems'' \cite[]{cvxguide}. This statement refers to two cases which are important to distinguish:
\begin{itemize}
	\item An appropriate solver is available for the problem, but the time to build the model in memory and pass it to the solver is a bottleneck. In this case, users are directed to use a low-level interface to the solver in place of the AML.
	\item The problem specified is simply too difficult for available solvers, whether in terms of memory use or computation time. In this case, users are directed to consider reformulating the problem or implementing specialized algorithms for it.
\end{itemize}
Our focus in this section is on the \textit{first} case, which is somehow typical of the large-scale performance of open-source AMLs implemented in high-level languages, as will be demonstrated in the numerical experiments in this section.

This performance gap between commercial and open-source AMLs can be partially explained by considering the possibly very different motivations of their respective authors; however, we posit that there is a more technical reason. In languages such as MATLAB and Python there is no programmatic access to the language's highly optimized expression parser. Instead, to handle nonlinear expressions such as \texttt{y*sin(x)}, one must either overload both the multiplication operator and the \texttt{sin} function, which leads to the expensive creation of many temporary objects as previously discussed, or manually parse the expression as a string, which itself may be slow and breaks the connection with the surrounding language. YALMIP and Pyomo implement the operator overloading approach, while CVX implements the string-parsing approach.

Julia, on the other hand, provides first-class access to its expression parser through its previously discussed metaprogramming features, which facilitates the generation of resulting code with performance comparable to that of commercial AMLs.  In Section~\ref{sec:nlpimpl} we describe our proof-of-concept implementation, followed by computational results in Section~\ref{sec:nlpresults}.

\subsection{Implementation in Julia}\label{sec:nlpimpl}

Whereas linear expressions are represented as sparse vectors of nonzero values, nonlinear expressions are represented as algebraic \textit{expression graphs}, as will be later illustrated. Expression graphs, when available, are integral to the solution of nonlinear optimization problems. The AML is responsible for using these graphs to evaluate function values and first and second derivatives as requested by the solver (typically through \textit{callbacks}). Additionally, they may be used by the AML to infer problem structure in order to decide which solution methods are appropriate~\cite[]{DrAmpl} or by the solver itself to perform important problem reductions in the case of mixed-integer nonlinear programming~\cite[]{couenne}.

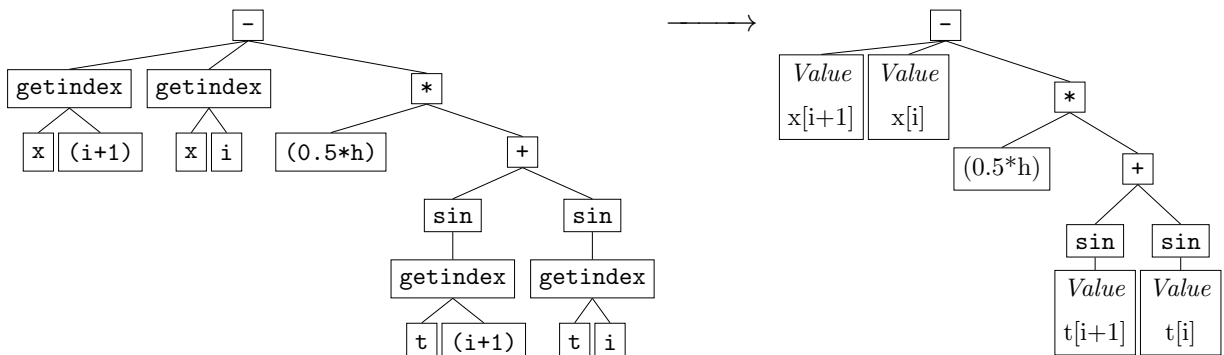
\begin{figure}[t]
	\centering
	\raisebox{-\height}{
\begin{tikzpicture}[scale=0.8]
	\tikzset{every internal node/.style={rectangle,draw}}
	\tikzset{every tree node/.style={font=\tt}}
	\tikzset{every leaf node/.style={rectangle,draw,minimum height=18pt}}
	\Tree [.- [.getindex x (i+1) ] [.getindex x i ] [.* (0.5*h) [.+ [.sin [.getindex t (i+1) ] ] [.sin [.getindex t {i} ] ] ] ] ]
\end{tikzpicture}}
\raisebox{-\height}{$\xrightarrow{\hspace*{1cm}}$}
\raisebox{-\height}{
\begin{tikzpicture}[scale=0.8]
	\tikzset{every internal node/.style={rectangle,draw,font=\tt}}
	\tikzset{every tree node/.style={align=center}}
	\tikzset{every leaf node/.style={rectangle,draw}}
	\tikzset{level distance=34pt}
	\Tree [.- {\textit{Value}\\x[i+1]} {\textit{Value}\\x[i]} [.* (0.5*h) [.+ [.sin {\textit{Value}\\t[i+1]} ] [.sin {\textit{Value}\\t[i]} ] ] ] ] 
\end{tikzpicture}}

\caption{A macro called with the text expression \texttt{x[i+1] - x[i] - (0.5h)*(sin(t[i+1])+sin(t[i]))} is given as input as the expression tree on the left. Parentheses indicate subtrees combined for brevity. At this stage, \texttt{symbols} are abstract and not resolved. To prepare the nonlinear expression, the macro produces code that generates the expression tree on the right with variable placeholders spliced in from the runtime context. }\label{fig:nlmacro}
\end{figure}

Analogously to the linear case, where macros are used to generate code which forms sparse vector representations, a macro was implemented which generates code to form nonlinear expression trees. Macros, when called, are provided an expression tree of the input; however, symbols are not resolved to values. Indeed, values do not exist at compile time when macros are evaluated. The task of the macro, therefore, is to generate code which replicates the input expression with runtime values (both numeric constants and variable placeholder objects) spliced in, as illustrated in Figure~\ref{fig:nlmacro}. This splicing of values is by construction and does not require expensive runtime calls such as MATLAB's \texttt{eval} function.

The implementation is compact, approximately 20 lines of code including support for the \texttt{sum\{\}} syntax presented in Section~\ref{sec:mathprogdesign}. While a nontrivial understanding of Julia's metaprogramming syntax is required to implement such a macro, the effort should be compared with what would be necessary to obtain the equivalent output and performance from a low-level language; in particular, one would need to write a custom expression parser.

Given expression trees for the constraints~\eqref{eq:constr}, we consider computing the Jacobian matrix 
\[
J(x) = \left[\begin{array}{c} \nabla g_1(x) \\ \nabla g_2(x) \\ \vdots \\ \nabla g_m(x)\end{array}\right],
\]
where $\nabla g_i(x)$ is a row-oriented gradient vector.
Unlike the typical approach of using automatic differentiation for computing derivatives in AMLs~\cite[]{Gay96moreAD}, a simpler method based on symbolic differentiation can be equally as efficient in Julia. In particular, we derive the form of the \textit{sparse} Jacobian matrix by applying the chain rule symbolically and then, using JIT compilation, \textit{compile a function which evaluates the Jacobian} for any given input vector. This process is accelerated by identifying equivalent expression trees (those which are symbolically identical and for which there exists a one-to-one correspondence between the variables present) and only performing symbolic differentiation once per equivalent expression.

The implementation of the Jacobian computation spans approximately 250 lines of code, including the basic logic for the chain rule. In the following section it is demonstrated that evaluating the Jacobian using the JIT compiled function is as fast as using AMPL through the low-level amplsolver library~\cite[]{amplsolver}, presently a de-facto standard for evaluating derivatives in nonlinear models. Interestingly, there is an executable accompanying the amplsolver library (\texttt{nlc}) which generates and compiles C code to evaluate derivatives for a specific model, although it is seldom used in practice because of the cost of compilation and the marginal gains in performance. However, in a language such as Julia with JIT compilation, compiling functions generated at runtime can be a technique to both simplify an implementation and obtain performance comparable to that of low-level languages. 

\subsection{Computational tests}\label{sec:nlpresults}

We test our implementation on two nonlinear optimization problems obtained from Hans Mittelmann's AMPL-NLP benchmark set ({\small\url{http://plato.asu.edu/ftp/ampl-nlp.html}}). Experiments were performed on a Linux system with an Intel Xeon E5-2650 processor. Note that we have not developed a complete nonlinear AML; the implementation is intended to serve as a proof of concept only. The operations considered are solely the construction of the model and the evaluation of the Jacobian of the constraints. Hence, objective functions and right-hand side expressions are omitted or simplified below.

The first instance is \textit{clnlbeam}:
\begin{align*}
	\min_{t,x,u \in \mathbb{R}^{n+1}} \quad & ...\\
	\text{subject to}\quad& x_{i+1} - x_{i} - \frac{1}{2n}(\sin(t_{i+1}) + \sin(t_i)) = 0\quad i = 1,\ldots,n\\
		& t_{i+1} - t_i - \frac{1}{2n}u_{i+1} - \frac{1}{2n}u_i = 0\quad i=1,\ldots,n\\
	 &-1 \leq t_i \leq 1,\quad -0.05 \leq x_i \leq 0.05\quad i = 1,\ldots,n+1
\end{align*}

We take $n =$ 5,000, 50,000, and 500,000. The following code builds the corresponding model in Julia using our proof-of-concept implementation:

\vspace{0.1in}
\begin{MyCode}
m = Model(:Min)
h = 1/n
@defVar(m, -1 <= t[1:(n+1)] <= 1)
@defVar(m, -0.05 <= x[1:(n+1)] <= 0.05)
@defVar(m, u[1:(n+1)])

for i in 1:n
	@addNLConstr(m, x[i+1] - x[i] - 
		(0.5h)*(sin(t[i+1])+sin(t[i])) == 0)
end
for i in 1:n
	@addNLConstr(m, t[i+1] - t[i] - 
		(0.5h)*u[i+1] - (0.5h)*u[i] == 0)
end
\end{MyCode}

The second instance is \textit{cont5\_1}:
\begin{align*}
	\min_{y \in \mathbb{R}^{(n+1)\times (n+1)},u \in \mathbb{R}^n} \quad & ...\\
	\text{subject to}\quad& n(y_{i+1,j+1} - y_{i,j+1}) - a(y_{i,j} - 2y_{i,j} + y_{i,j+1} + y_{i+1,j} - 2y_{i+1,j+1} + y_{i+1,j+2}) = 0\\
						  & \quad i=1,\ldots,n, \, j = 1,\ldots,n-1\\
						  & y_{i+1,3} - 4y_{i+1,2} + 3y_{i+1,1} = 0\quad i=1,\ldots,n\\
					   & c(y_{i+1,n-1} - 4y_{i+1,n} + 3y_{i+1,n+1}) + y_{i+1,n+1} - u_i + y_{i+1,n+1}((y_{i+1,n+1})^2)^\frac{3}{2} = 0 \\
						& \quad i = 1,\ldots,n,
\end{align*}

where $a = \frac{8n^2}{\pi^2}$ and $c = \frac{2n}{\pi}$. We take $n =$ 200, 400, and 1,000.

The dimensions of these instances and the number of nonzero elements in the corresponding Jacobian matrices are listed in Table~\ref{tab:nldims}. In Table~\ref{tab:nlresults} we present a benchmark of our implementation compared with AMPL, YALMIP (MATLAB), and Pyomo (Python). JIT compilation of the Jacobian function is included in the 	``Build model'' phase for Julia. Observe that Julia performs as fast as AMPL, if not faster. Julia's advantage over AMPL is partly explained by AMPL's need to write the model to an intermediate \texttt{nl} file before evaluating Jacobians; this I/O time is included. AMPL's preprocessing features are disabled. YALMIP performs well on the mostly linear cont5\_1 instance but is unable to process the largest clnlbeam instance in under an hour. Pyomo's performance is more consistent but over 50x slower than Julia on the largest instances. Pyomo is run under pure Python; it does not support JIT accelerators such as PyPy. 

\begin{table}[ht]
	\centering
	\caption{Nonlinear test instance dimensions. Nz = Nonzero elements in Jacobian matrix.}
	\label{tab:nldims}
	\begin{tabular}{lrrr}
		\toprule
		Instance & \# Vars. & \# Constr. & \# Nz\\
		\cmidrule(l){1-1}\cmidrule(l){2-2}\cmidrule(l){3-3}\cmidrule(l){4-4}
		clnlbeam-5 & 15,003 & 10,000 & 40,000\\
		clnlbeam-50 & 150,003 & 100,000 & 400,000\\
		clnlbeam-500 & 1,500,003 & 1,000,000 & 4,000,000\\
		cont5\_1-2 & 40,601 & 40,200 & 240,200 \\
		cont5\_1-4 & 161,201 & 160,400 & 960,400\\
		cont5\_1-10 & 1,003,001 & 1,001,000 & 6,001,000\\
		\bottomrule\\
	\end{tabular}
\end{table}

\begin{table}[ht]
	\centering
	\caption{Nonlinear benchmark results. ``Build model`` includes writing and reading model files, if required, and precomputing the structure of the Jacobian. Pyomo uses AMPL for Jacobian evaluations. }
	\label{tab:nlresults}
	\begin{tabular}{lrrrrrrr}
		\toprule
		&\multicolumn{4}{c}{Build model (s)} & \multicolumn{3}{c}{Evaluate Jacobian (ms)} \\
		\cmidrule(l){2-5}\cmidrule(l){6-8}
		Instance & AMPL & Julia & YALMIP & Pyomo & AMPL & Julia & YALMIP \\
		\cmidrule(l){1-1}\cmidrule(l){2-2}\cmidrule(l){3-3}\cmidrule(l){4-4}\cmidrule(l){5-5}\cmidrule(l){6-6}\cmidrule(l){7-7}\cmidrule(l){8-8}
		clnlbeam-5 & 0.2 & 0.1	& 36.0 & 2.3 & 0.4& 0.3 & 8.3\\
		clnlbeam-50 & 1.8 & 0.3	& 1344.8& 23.7 & 7.3 & 4.2 & 96.4\\
		clnlbeam-500 & 18.3 & 3.3 & \textgreater 3600 & 233.9 & 74.1 & 74.6 & *\\
		cont5\_1-2 & 1.1 &  0.3 & 2.0 & 12.2	& 1.1 &  0.8 & 9.3\\
		cont5\_1-4 & 4.4  & 1.4 & 1.9 & 49.4	& 5.4  & 3.0 & 37.4\\
		cont5\_1-10 & 27.6 & 6.1 & 13.5	& 310.4 & 33.7 & 39.4 & 260.0\\
		\bottomrule
		\\
	\end{tabular}
\end{table}

%% file: simplex.tex
\section{Implementing Optimization Algorithms}
\label{sec:Simplex}

In this section we evaluate the performance of Julia for implementation of the simplex method for linear programming, arguably one of the most important algorithms in the field of operations research. Our aim is not to develop a complete implementation but instead to compare the performance of Julia to that of other popular languages, both high- and low-level, on a benchmark of core operations. 

Although high-level languages can achieve good performance when performing \textit{vectorized} operations (that is, block operations on dense vectors and matrices), state-of-the-art implementations of the simplex method are characterized by their effective exploitation of sparsity (the presence of many zeros) in all operations, and hence, they use sparse linear algebra. Opportunities for vectorized operations are small in scale and do not represent a majority of the execution time; see~\cite{HallParSimplex}. Furthermore, the sparse linear algebra operations used, such as~\cite{SuhlLU}'s $LU$ factorization, are specialized and not provided by standard libraries.

The simplex method is therefore an example of an algorithm that requires a low-level coding style, in particular, manually-coded loops, which are known to have poor performance in languages such as Matlab or Python (see, e.g., ~\cite{numpy}). To achieve performance in such cases, one would be required to code time-consuming loops in another language and link to these separate routines from the high-level language, using, for example, Matlab's MEX interface. Our benchmarks will demonstrate, however, that within Julia, the native performance of this style of computation can nearly achieve that of low-level languages.

\subsection{Benchmark Operations}
A presentation of the simplex algorithm and a discussion of its computational components are beyond the scope of this paper. We refer the reader to~\cite{MarosBook} and~\cite{KobersteinThesis} for a comprehensive treatment of modern implementations, which include significant advances over versions presented in most textbooks. We instead present three selected operations from the \textit{revised dual simplex} method in a mostly self-contained manner. Knowledge of the simplex algorithm is helpful but not required. The descriptions are realistic and reflect the essence of the routines as they might be implemented in an efficient implementation. 

The first operation considered is a matrix-transpose-vector product (Mat-Vec). In the revised simplex method, this operation is required in order to form a row of the \textit{tableau}. A nonstandard aspect of this Mat-Vec is that we would like to consider the matrix formed by a constantly changing subset the columns (those corresponding to the non-basic variables). Another important aspect is the treatment of sparsity of the \textit{vector} itself, in addition to that of the matrix~\cite[]{HallHyperSparse}. This is achieved algorithmically by using the nonzero elements of the vector to form a linear combination of the rows of the matrix, instead of the more common approach of computing dot-products with the columns, as illustrated in \eqref{eq:matvec}. This follows from viewing the matrix $A$ equivalently as either a collection of column vectors $A_i$ and or as row vectors $a_i^T$.
\begin{equation}\label{eq:matvec}
A = \left[\begin{array}{cccc} A_1 & A_2 & \cdots & A_n\end{array}\right] = \left[\begin{array}{c} a_1^T \\ a_2^T \\ \vdots \\ a_m^T\end{array}\right]\longrightarrow 
A^Tx = \left[\begin{array}{c} A_1^Tx \\ A_2^Tx \\ \vdots \\ A_n^Tx \end{array}\right] = \sum_{\substack{i=1 \\ x_i \neq 0}}^m a_ix_i
\end{equation}

\begin{algorithm}[ht]\small
\caption{Restricted sparse matrix transpose-dense vector product}\label{alg:matvec}
\begin{algorithmic}
\State \textbf{Input:} Sparse column-oriented $m \times n$ matrix $A$, dense vector $x \in \mathbb{R}^m$, and \\
\hspace{1.34cm} flag vector $\mathcal{N} \in \{0,1\}^n$ (with $n-m$ nonzero elements)
\State \textbf{Output:} $y := A_{\mathcal{N}}^Tx$ as a dense vector, where $\mathcal{N}$ selects columns of $A$
\For{$i$ \textbf{in} $\{1,\ldots,n\}$}
\If{$\mathcal{N}_i = 1$} 
\State $s \leftarrow 0$
\ForAll{nonzero element $q$ (in row $j$) of $i$th column of $A$}
\State $s \leftarrow s + q*x_j$ \Comment{Compute dot-product of $x$ with column $i$}
\EndFor
\State $y_i \leftarrow s$
\EndIf
\EndFor
\end{algorithmic}
\end{algorithm}

\begin{algorithm}[ht]\small
\caption{Sparse matrix transpose-sparse vector product}\label{alg:matvec2}
\begin{algorithmic}
\State \textbf{Input:} Sparse row-oriented $m \times n$ matrix $A$ and sparse vector $x \in \mathbb{R}^m$
\State \textbf{Output:} Sparse representation of $A^Tx$.
\ForAll{nonzero element $p$ (in index $j$) in $x$}
\ForAll{nonzero element $q$ (in column $i$) of $j$th row of $A$}
\State Add $p*q$ to index $i$ of output. \Comment{Compute linear combination of rows of $A$} 
\EndFor
\EndFor
\end{algorithmic}
\end{algorithm}

The Mat-Vec operation is illustrated for dense vectors in Algorithm~\ref{alg:matvec} and for sparse vectors in Algorithm~\ref{alg:matvec2}. Sparse matrices are provided in either compressed sparse column (CSC) or compressed sparse row (CSR) format as appropriate~\cite[]{DuffCSC}. Note that in Algorithm~\ref{alg:matvec} we use a flag vector to indicate the selected columns of the matrix $A$. This corresponds to skipping particular dot products. The result vector has a memory layout with $n$, not $n-m$ entries. This form could be desired in some cases for subsequent operations and is illustrative of the common practice in simplex implementations of designing data structures with a global view of the operations in which they will be used~\cite[Chap. 5]{MarosBook}. In Algorithm~\ref{alg:matvec2} we omit what would be a costly flag check for each nonzero element of the row-wise matrix; the gains of exploiting sparsity often outweigh the extra floating-point operations. 
\begin{algorithm}[ht]\small
\caption{Two-pass stabilized minimum ratio test (dual simplex)}\label{alg:ratio}
\begin{multicols}{2}
\begin{algorithmic}
\State \textbf{Input:} Vectors $d, \alpha \in \mathbb{R}^n$, state vector $s \in \{$``lower'',``basic''$\}^n$, parameters $\epsilon_{\text{P}}, \epsilon_{\text{D}} > 0$
\State \textbf{Output:} Solution index $result$.
\State $\Theta_{\max} \leftarrow \infty$
\For{$i$ \textbf{in} $\{1,\ldots,n\}$}
\If{$s_i =$ ``lower'' and $\alpha_i > \epsilon_{\text{P}}$}
\State Add index $i$ to list of candidates
\State $\Theta_{\max} \leftarrow \min(\frac{d_i+\epsilon_{\text{D}}}{\alpha_i},\Theta_{\max})$
\EndIf
\EndFor
\State $\alpha_{\max} \leftarrow 0$, $result \leftarrow 0$
\For{$i$ \textbf{in} list of candidates}
\If {$d_i/\alpha_i \leq \Theta_{\max}$ and $\alpha_i > \alpha_{\max}$}
\State $\alpha_{\max} \leftarrow \alpha_i$
\State $result \leftarrow i$
\EndIf
\EndFor
\end{algorithmic}
\end{multicols}
\end{algorithm}

The second operation is the minimum ratio test, which determines both the step size of the next iteration and the constraint that prevents further progress. Mathematically this may be expressed as
\[
	\min_{\alpha_i > 0} \frac{d_i}{\alpha_i},
\]
for given vectors $d$ and $\alpha$. While seemingly simple, this operation is one of the more complex parts of an implementation, as John Forrest mentions in a comment in the source code of the open-source Clp solver. We implement a relatively simple two-pass variant (Algorithm~\ref{alg:ratio}) due to~\cite{HarrisDevex} and described more recently in~\cite[Sect. 6.2.2.2]{KobersteinThesis}, whose aim is to avoid numerical instability caused by small values of $\alpha_i$. In the process, small infeasibilities up to a numerical tolerance $\epsilon_{\text{D}}$ may be created. Note that our implementation handles both upper and lower bounds; Algorithm~\ref{alg:ratio} is simplified in this respect for brevity. A sparse variant is easily obtained by looping over the nonzero elements of $\alpha$ in the first pass.

The third operation is a modified form of the vector update $y \leftarrow \alpha x + y$ (Axpy). In the variant used in the simplex algorithm, the value of each updated component is tested for membership in an interval. For example, given a tolerance $\epsilon$, a component belonging to the interval $(-\infty,-\epsilon)$ may indicate loss of numerical feasibility, in which case a certain corrective action, such as local perturbation of problem data, may be triggered. This procedure is more naturally expressed using an explicit loop over elements of $x$ instead of performing operations on vectors.

The three operations discussed represent a nontrivial proportion of execution time of the simplex method, between 20\% and 50\% depending on the problem instance~\cite[]{HallHyperSparse}. Most of the remaining execution time is spent in factorizing and solving linear systems using specialized procedures, which we do not implement because of their complexity. 
\subsection{Results}

The benchmark operations described in the previous section were implemented in Julia, C++, MATLAB, and Python. Examples of code have been omitted for brevity. The style of the code in Julia is qualitatively similar to that of the other high-level languages. Readers are encouraged to view the implementations available in the online supplement. To measure the overhead of bounds-checking, a validity check performed on array indices in high-level languages, we implemented a variant in C++ with explicit bounds checking. We also consider executing the Python code under the PyPy engine~\cite[]{PyPy}, a JIT-compiled implementation of Python. We have not used the popular NumPy library in Python because it would not alleviate the need for manually coded loops and so would provide little speed benefit. No special runtime parameters are used, and the C++ code is compiled with \texttt{-O2}.

Realistic input data were generated by running a modified implementation of the dual simplex algorithm on a small set of standard LP problems and recording the required input data for each operation from iterations sampled uniformly over the course of the algorithm. At least 200 iterations are recorded from each instance. Using such data from real instances is important because execution times depend significantly on the sparsity patterns of the input. The instances we consider are greenbea, stocfor3, and ken-13 from the NETLIB repository~\cite[]{gay1985electronic} and the fome12 instance from Hans Mittelmann's benchmark set~\cite[]{mittelmann2009benchmarks}. These instances represent a range of problem sizes and sparsity structures. 

Experiments were performed under the Linux operating system on a laptop with an Intel i5-3320M processor.  See Table~\ref{tab:simplex} for a summary of results. Julia consistently performs within a factor of 2 of the implementation in C++ with bounds checking, while MATLAB and PyPy are within a factor of 4 to 18. Pure Python is far from competitive, being at least 70x slower than C++.

\begin{figure}[ht]
	\centering
	\includegraphics{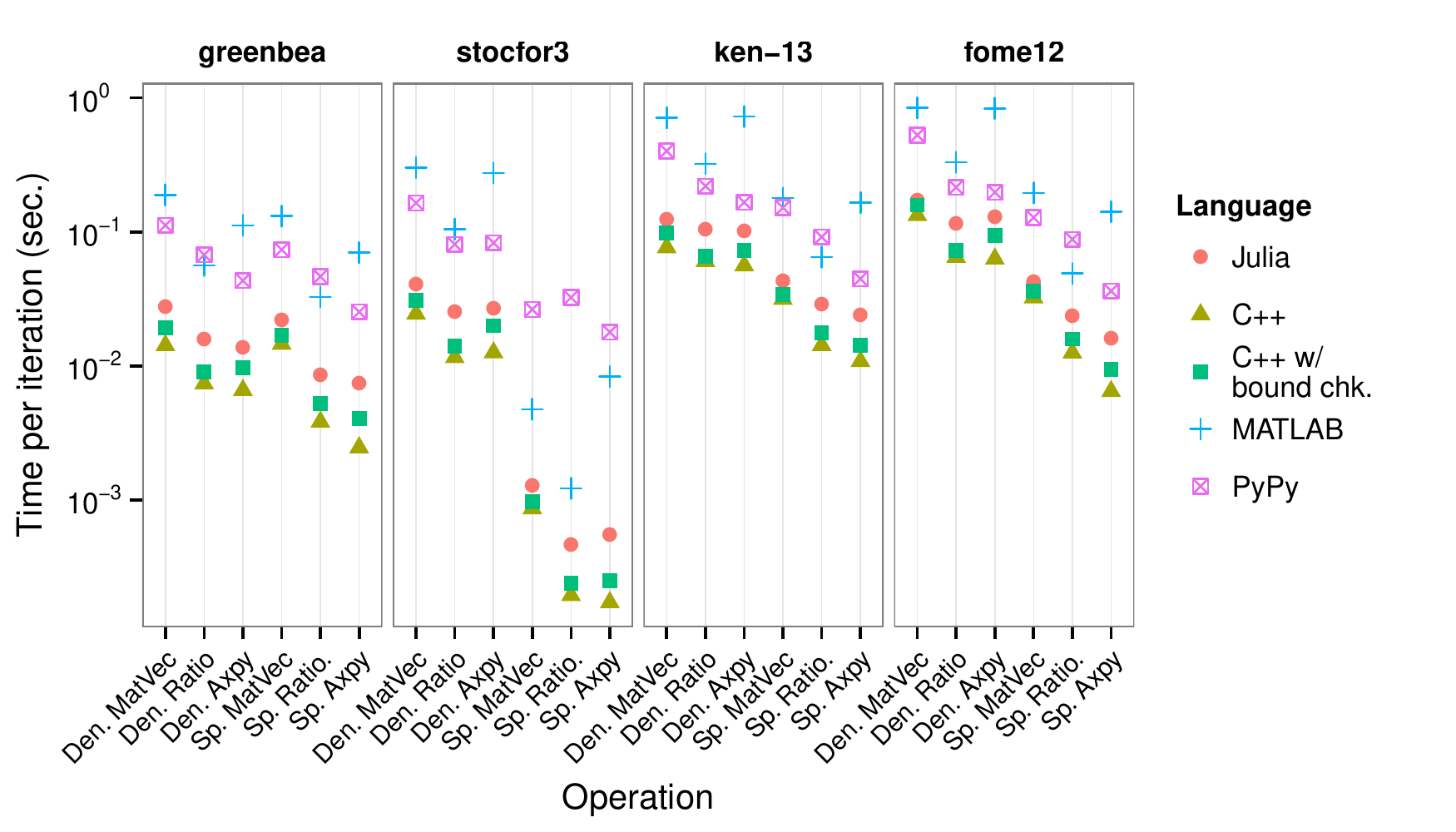}
	\caption{Average execution time for each operation and language, by instance. Compared with MATLAB and PyPy, the execution time of Julia is significantly closer to that of C++.}\label{fig:simplex}
\end{figure}

Figure~\ref{fig:simplex} displays the absolute execution times broken down by instance. We observe the consistent performance of Julia, while that of MATLAB and PyPy are subject to more variability. In all cases except the smaller greenbea instance, use of the vector-sparse routines significantly decreases execution time, although PyPy's performance is relatively poorer on these routines.

\begin{table}[th]
	\centering
	\caption{Execution time of each language (version listed below) relative to C++ with bounds checking. Lower values are better. Figures are geometric means of average execution times over iterations over 4 standard LP problems. Recorded value is fastest time of three repetitions. Dense/sparse distinction refers to the \textit{vector} $x$; all matrices are sparse. }\label{tab:simplex} 
	\begin{tabular}{llrrrrr}
		\toprule
		&& Julia & C++ & MATLAB & PyPy & Python \\
		&Operation& 0.1 & GCC 4.7.2 & R2012b & 1.9 & 2.7.3\\
		\cmidrule(l){2-2}\cmidrule(l){3-3}\cmidrule(l){4-4}\cmidrule(l){5-5}\cmidrule(l){6-6}\cmidrule(l){7-7}
		Dense &		Mat-Vec $(A^T_{\mathcal{N}}x)$& 1.27&	0.79&	7.78&	4.53&	84.69\\	
		
			   &Min. ratio test&	1.67&	0.86&	5.68&	4.54&	70.95\\
			   &Axpy ($y \leftarrow \alpha x + y$) &	1.37	&0.68	&10.88	&3.07	&83.71\medskip\\	
		Sparse &Mat-Vec ($A^Tx$)&	1.25&	0.89&	5.72&	6.56&	69.43	\\
				&Min. ratio test &	1.65&	0.78&	4.35&	13.62&	73.47	\\
		  &Axpy ($y \leftarrow \alpha x + y$) & 1.84&	0.68&	17.83&	8.57&	81.48\\
		\bottomrule
	\end{tabular}
\end{table}

Our results are qualitatively similar to those reported by~\cite{JuliaArxiv} on a set of unrelated general language benchmarks and thus serve as an independent corroboration of their findings that Julia's performance is within a factor of 2 of equivalent low-level compiled code.

%% file: main.bbl
\begin{thebibliography}{24}
\expandafter\ifx\csname natexlab\endcsname\relax\def\natexlab#1{#1}\fi
\expandafter\ifx\csname url\endcsname\relax
  \def\url#1{{\tt #1}}\fi
\expandafter\ifx\csname urlprefix\endcsname\relax\def\urlprefix{URL }\fi
\expandafter\ifx\csname urlstyle\endcsname\relax
  \expandafter\ifx\csname doi\endcsname\relax
  \def\doi#1{doi:\discretionary{}{}{}#1}\fi \else
  \expandafter\ifx\csname doi\endcsname\relax
  \def\doi{doi:\discretionary{}{}{}\begingroup \urlstyle{rm}\Url}\fi \fi

\bibitem[{Belotti et~al.(2009)Belotti, Lee, Liberti, Margot, and
  W\"{a}chter}]{couenne}
Belotti, Pietro, Jon Lee, Leo Liberti, François Margot, Andreas W\"{a}chter.
  2009.
\newblock Branching and bounds tightening techniques for non-convex {MINLP}.
\newblock {\it Optimization Methods and Software\/} {\bf 24} 597--634.

\bibitem[{Bezanson et~al.(2012)Bezanson, Karpinski, Shah, and
  Edelman}]{JuliaArxiv}
Bezanson, Jeff, Stefan Karpinski, Viral~B. Shah, Alan Edelman. 2012.
\newblock Julia: A fast dynamic language for technical computing.
\newblock {\it CoRR\/} {\bf abs/1209.5145}.

\bibitem[{Bixby(2002)}]{bixby2002solving}
Bixby, Robert~E. 2002.
\newblock Solving real-world linear programs: A decade and more of progress.
\newblock {\it Operations research\/} {\bf 50} 3--15.

\bibitem[{Bolz et~al.(2009)Bolz, Cuni, Fijalkowski, and Rigo}]{PyPy}
Bolz, Carl~Friedrich, Antonio Cuni, Maciej Fijalkowski, Armin Rigo. 2009.
\newblock Tracing the meta-level: {PyPy}'s tracing {JIT} compiler.
\newblock {\it Proceedings of the 4th workshop on the Implementation,
  Compilation, Optimization of Object-Oriented Languages and Programming
  Systems\/}. ICOOOLPS '09, ACM, New York, 18--25.

\bibitem[{Brooke et~al.(1999)Brooke, Kendrick, Meeraus, and Raman}]{GAMS}
Brooke, A., D.~Kendrick, A.~Meeraus, R.~Raman. 1999.
\newblock {\it {GAMS}: A User's Guide\/}.
\newblock Scientific Press.

\bibitem[{Duff et~al.(1989)Duff, Grimes, and Lewis}]{DuffCSC}
Duff, I.~S., Roger~G. Grimes, John~G. Lewis. 1989.
\newblock Sparse matrix test problems.
\newblock {\it ACM Trans. Math. Softw.\/} {\bf 15} 1--14.

\bibitem[{Fourer and Orban(2010)}]{DrAmpl}
Fourer, R., Dominique Orban. 2010.
\newblock {DrAmpl}: a meta solver for optimization problem analysis.
\newblock {\it Computational Management Science\/} {\bf 7} 437--463.

\bibitem[{Fourer et~al.(1993)Fourer, Gay, and Kernighan}]{AMPLBook}
Fourer, Robert, David~M Gay, Brian~W Kernighan. 1993.
\newblock {\it AMPL\/}.
\newblock Scientific Press.

\bibitem[{Gay(1985)}]{gay1985electronic}
Gay, David~M. 1985.
\newblock Electronic mail distribution of linear programming test problems.
\newblock {\it Mathematical Programming Society COAL Newsletter\/} {\bf 13}
  10--12.

\bibitem[{Gay(1996)}]{Gay96moreAD}
Gay, David~M. 1996.
\newblock More {AD} of nonlinear {AMPL} models: Computing hessian information
  and exploiting partial separability.
\newblock {\it in Computational Differentiation: Applications, Techniques, and
  Tools\/}. SIAM, 173--184.

\bibitem[{Gay(1997)}]{amplsolver}
Gay, David~M. 1997.
\newblock Hooking your solver to {AMPL}.
\newblock Tech. rep., Bell Laboratories, Murray Hill, NJ.

\bibitem[{Grant and Boyd(2013)}]{cvxguide}
Grant, Michael~C., Stephen~P. Boyd. 2013.
\newblock The {CVX} users' guide (release 2.0).
\newblock \urlprefix\url{http://cvxr.com/cvx/doc/CVX.pdf}.

\bibitem[{Hall(2010)}]{HallParSimplex}
Hall, J. 2010.
\newblock Towards a practical parallelisation of the simplex method.
\newblock {\it Computational Management Science\/} {\bf 7} 139--170.

\bibitem[{Hall and McKinnon(2005)}]{HallHyperSparse}
Hall, J., K.~McKinnon. 2005.
\newblock Hyper-sparsity in the revised simplex method and how to exploit it.
\newblock {\it Computational Optimization and Applications\/} {\bf 32}
  259--283.

\bibitem[{Harris(1973)}]{HarrisDevex}
Harris, Paula M.~J. 1973.
\newblock Pivot selection methods of the {DEVEX} {LP} code.
\newblock {\it Mathematical Programming\/} {\bf 5} 1--28.

\bibitem[{Hart et~al.(2011)Hart, Watson, and Woodruff}]{Pyomo}
Hart, William~E, Jean-Paul Watson, David~L Woodruff. 2011.
\newblock Pyomo: modeling and solving mathematical programs in {Python}.
\newblock {\it Mathematical Programming Computation\/} {\bf 3} 219--260.

\bibitem[{Koberstein(2005)}]{KobersteinThesis}
Koberstein, Achim. 2005.
\newblock The dual simplex method, techniques for a fast and stable
  implementation.
\newblock Ph.D. thesis, Universit\"{a}t Paderborn, Paderborn, Germany.

\bibitem[{Lattner and Adve(2004)}]{lattner2004llvm}
Lattner, Chris, Vikram Adve. 2004.
\newblock {LLVM}: A compilation framework for lifelong program analysis \&
  transformation.
\newblock {\it Code Generation and Optimization, 2004. International Symposium
  on\/}. IEEE, 75--86.

\bibitem[{Lofberg(2004)}]{YALMIP}
Lofberg, John. 2004.
\newblock {YALMIP}: A toolbox for modeling and optimization in {MATLAB}.
\newblock {\it Computer Aided Control Systems Design, 2004 IEEE International
  Symposium on\/}. IEEE, 284--289.

\bibitem[{Maros(2003)}]{MarosBook}
Maros, Istv\'{a}n. 2003.
\newblock {\it Computational Techniques of the Simplex Method\/}.
\newblock Kluwer Academic Publishers, Norwell, MA.

\bibitem[{Mitchell et~al.(2011)Mitchell, O'Sullivan, and Dunning}]{PULP}
Mitchell, Stuart, Michael O'Sullivan, Iain Dunning. 2011.
\newblock Pulp: A linear programming toolkit for python
  \urlprefix\url{https://code.google.com/p/pulp-or/}.
\newblock Unpublished manuscript.

\bibitem[{Mittelmann(2013)}]{mittelmann2009benchmarks}
Mittelmann, Hans. 2013.
\newblock Benchmarks for optimization software.
\newblock \urlprefix\url{http://plato.la.asu.edu/bench.html}.
\newblock Accessed April 28, 2013.

\bibitem[{Suhl and Suhl(1990)}]{SuhlLU}
Suhl, Uwe~H., Leena~M. Suhl. 1990.
\newblock Computing sparse {LU} factorizations for large-scale linear
  programming bases.
\newblock {\it ORSA Journal on Computing\/} {\bf 2} 325.

\bibitem[{van~der Walt et~al.(2011)van~der Walt, Colbert, and
  Varoquaux}]{numpy}
van~der Walt, S., S.C. Colbert, G.~Varoquaux. 2011.
\newblock The {NumPy} array: A structure for efficient numerical computation.
\newblock {\it Computing in Science Engineering\/} {\bf 13} 22 --30.

\end{thebibliography}
